\documentclass[a4paper,11pt, twoside]{article}
\usepackage[latin1]{inputenc}
\usepackage{amsfonts}
\usepackage{amsmath}
\usepackage{enumerate}
\usepackage{amssymb}
\usepackage{mdwlist}

\begin{document}
\title{On the Riemann Hypothesis and its generalizations}
\author{Daniel E. Borrajo Gutiérrez}
\date{\today}
\maketitle
\begin{abstract}
A proof for the original Riemann hypothesis is proposed based on the infinite Hadamard product representation for the Riemann zeta function and later generalized to Dirichlet L-functions. The extension of the hypothesis to other functions is also discussed.
\end{abstract}
\section{Introduction}
The Riemann zeta function $\zeta(s)$ is defined for $\operatorname{Re}(s) > 1$ as the infinite sum

\begin{equation}\label{zeta}
\zeta(s) = \sum_{n=1}^{\infty} {\frac1{n^s}},
\end{equation}
and can be extended by analytic continuation to any complex number $s \neq 1$, where \eqref{zeta} does not converge. It satisfies the following functional equation \cite{1} (p. 169):

\begin{equation}\label{functional}
\zeta(s) = 2^s \pi^{s-1} \sin{\left( \frac{\pi s}{2} \right)} \Gamma{(1-s)} \zeta{(1-s)}.
\end{equation}

From this equation it is easily seen that the function has infinite zeros where the sine function vanishes, i.e., negative even integers:

\begin{equation*}
\zeta{(-2l)} = 0 \quad \forall \ l \in \mathbb{N},
\end{equation*}

while positive integers cannot be considered in \eqref{functional} because of the poles of the Gamma function. These are called \emph{trivial zeros} and are of little interest in this paper. The interesting ones are the non-trivial zeros, which lie in the \emph{critical strip} $0 < \operatorname{Re}(s) < 1$ (from \eqref{zeta} and \eqref{functional} it can be seen that non-trivial zeros cannot be found outside the critical strip.)

Riemann himself obtained 3 of these zeros with $\operatorname{Re}(\rho) = 1/2$ and proposed in 1859 the following:

\begin{center}
\emph{The non-trivial zeros of the zeta function have all real part 1/2 (RH)}.
\end{center}

Since then, many proofs and disproofs of this hypothesis have been proposed, but so far none has been widely accepted as correct. In the next section, a simple proof is proposed. We will work only on the critical strip, and call the non-trivial zeros simply 'zeros'.

\section{Proposed proof}
We start from the symmetrical functional equation
\begin{equation}\label{sym}
\Gamma{\left( \frac{s}{2} \right)} \pi^{-s/2} \zeta{(s)} = \Gamma{\left( \frac{1-s}{2} \right)} \pi^{-(1-s)/2} \zeta{(1-s)},
\end{equation}
and define the function $ \Phi{(s)} $ as the left-hand side of \eqref{sym}, which satisfies $\Phi{(s)} = \Phi{(1-s)}$. We further define
\begin{equation}
F(t) \equiv \Phi{(1/2+it)},
\end{equation}
and assert that $F(t)$ is real-valued due to the following property: \\

\emph{Property 1}: let $f(z)$ be a complex-valued meromorphic function whose restriction to the real numbers is real-valued; then:
\begin{equation}
f(\overline{z}) = \overline{f(z)}.
\end{equation}
\emph{Proof}: this can be seen from the Laurent series of such a function:
\begin{equation*}
f(z) = \sum_{n=-\infty}^{\infty} {a_n (z-z_0)^n};
\end{equation*}
One must have $ a_n = \overline{a}_n $ for, otherwise, its restriction to the real numbers would not be real-valued, so:
\begin{equation*}
f(\overline{z}) = \sum_{n=-\infty}^{\infty} {a_n (\overline{z-z_0})^n} = \sum_{n=-\infty}^{\infty} {\overline{a}_n (\overline{z-z_0})^n} = \overline{f(z)}. \quad \blacksquare
\end{equation*}
$F(t)$ satisfies both requirements because it is a product of functions which indeed satisfy them. So
\begin{equation*}
\Phi{(1/2+it)} = \Phi{(1/2-it)} \Longrightarrow F(t) = \overline{F(t)}.
\end{equation*}
To representate the zeta function, we use the Hadamard product:

\begin{equation}\label{zetaHadamard}
\zeta{(s)} = \frac{\pi^{s/2}}{2(s-1)\Gamma \left( 1 + \frac{s}{2} \right)} \prod_\rho {\left( 1 - \frac{s}{\rho} \right)},
\end{equation}
where $\rho$ are the zeros of the zeta function. We obtain

\begin{equation}\label{f2}
F(t) = - \frac{B(t)}{\frac1{4} + t^2},
\end{equation}
where
\begin{equation}\label{zeros}
B(t) = \prod_\rho {\left( 1 - \frac{1/2+it}{\rho} \right)}
\end{equation}

and property $\Gamma(1+z) = z \Gamma(z)$ has been employed. We note that \eqref{f2} gives $B(t) \in \mathbb{R}$.

Let us suppose the Riemann hypothesis is not true. Because, according to property 1 $\zeta(\overline{s}) = \overline{\zeta(s)}$, the complex conjugate of a zero is also a zero, so the zeros outside the critical line should come in sets of four:

\begin{equation*}
\left\{ 
  \begin{array}{l l l l}
    \rho_1 = \sigma + i \tau, \\
    \rho_2 = 1-\sigma - i \tau, \\
    \rho_3 = \sigma - i \tau, \\
    \rho_4 = 1- \sigma + i \tau.
  \end{array} \right.
\end{equation*}

We now separate the Riemann zeros ($\operatorname{Re}(s)=1/2$) and express the non-Riemann ones in polar coordinates in \eqref{zeros}:

\begin{equation}\label{allzeros}
B(t) = \left[ \prod_{\tau} {\left( 1 - \frac{1/2 + it}{1/2 + i\tau} \right)} \right] \mathfrak{R} e^{i \phi} \in \mathbb{R}.
\end{equation}

The product over the Riemann zeros is real; in fact, if we combine the pairs of zeros $\rho$ and $1-\rho$, we obtain
\begin{equation} \label{Riemannzeros}
B_{RH}(t) = \prod_{\substack{ \sigma = 1/2 \\ \tau > 0}} { \frac{\tau^2 - t^2}{1/4 + \tau^2} } \in \mathbb{R}.
\end{equation}

Now, for each set, we perform the product over the four matching zeros and, because of \eqref{allzeros} and \eqref{Riemannzeros}:

\begin{equation*}
\phi(t) = \sum_{\substack{ \sigma < 1/2 \\ \tau > 0}} {\operatorname{atan} \left( \frac{2t(1-2\sigma)(t^2-\tau^2-\alpha)}{(t^2+\tau^2)^2+\alpha(\alpha - 3t^3 - \tau^2)} \right)} = n(t) \pi,
\end{equation*}

where $n(t) \in \mathbb{Z}$ and $\alpha = (1/2-\sigma)^2$. In fact, $n(t)$ takes a constant value for all $t$: suppose $n(t_1) = y_1$ and $n(t_2) = y_2$ and $y_1 \leq y_2$. By the intermediate value theorem, $n(t)$ must take all the values in $[y_1,y_2]$, but the range is the integers, so this is only true if $n$ only takes a single value. We also know that $\phi(0) = 0$ so, by virtue of continuity, this value is indeed 0:

\begin{equation}\label{polar}
\phi(t) = \sum_{\substack{ \sigma < 1/2 \\ \tau > 0}} {\operatorname{atan} \left( \frac{2t(1-2\sigma)(t^2-\tau^2-\alpha)}{(t^2+\tau^2)^2+\alpha(\alpha - 3t^3 - \tau^2)} \right)} = 0.
\end{equation}

To obtain \eqref{polar}, both numerator and denominator of the argument have been divided by the common denominator of the real and imaginary parts:

\begin{equation*}
D(\sigma, \tau) = \sigma^4 + \sigma^2(1-2\sigma+2\tau^2) + \tau^2(1-2\sigma) + \tau^4;
\end{equation*}

this is allowed provided $D(\sigma, \tau) \neq 0$ everywhere inside the critical strip. In fact, $D=0$ for $\sigma = \tau = 0$ and $D$ diverges for high values of $\tau$. In the middle, it can be shown that $D$ has no minimum value because

\begin{equation*}
\frac{\partial D(\sigma, \tau)}{\partial \sigma} = \frac{\partial D(\sigma, \tau)}{\partial \tau} = 0
\end{equation*}
is not satisfied for any real $\sigma$ and $\tau$.

Because $B(t) \in \mathbb{R} \ \forall \ t$, it suffices to find a $t$ for which \eqref{polar} does not hold. In fact, if $\rho_0 = \sigma_0 + i \tau_0$ is the zero (or set of zeros) with the lowest value of $\tau^2 + \alpha$, we can choose $t_0$ so that
\begin{equation*}
t_0^2 - \tau_0^2 - \left( \frac1{2} - \sigma_0 \right)^2 = 0;
\end{equation*}

for $B(t_0)$, the first term in \eqref{polar} vanishes and the subsequent terms are all negative because $t_0^2 < \tau^2 + (1/2-\sigma^2)$ and $\operatorname{atan}(x)$ is an even monotonic function. The only possibility is that $\rho_0$ be the only non-Riemann zero, but this is not allowed either because we can choose $t_1 > t_0$ so that the only term in \eqref{polar} does not vanish and we arrive at a contradiction.

The last step is to verify that the denominator $z(\sigma, \tau, t)$ in \eqref{polar} is not zero in the critical strip; in fact, $z(1/2,0,0) = 0$ and the condition

$$\frac{\partial z(\sigma, \tau, t)}{\partial \sigma} = \frac{\partial z(\sigma, \tau, t)}{\partial \tau} = \frac{\partial z(\sigma, \tau, t)}{\partial t} = 0$$

is not satisfied for any positive $\sigma$, $\tau$, $t$. This concludes the proof.  $\blacksquare$

\section{The Generalized Riemann Hypothesis}
In 1837, Peter Gustav Lejeune Dirichlet introduced the L-series, defined for $\operatorname{Re}(s) > 1$ as

\begin{equation}\label{DirichletL}
L(s, \chi) = \sum_{n=1}^{\infty} {\frac{\chi{(n)}}{n^s}},
\end{equation}
where $\chi(n), \ n \in \mathbb{Z}$ is the \emph{Dirichlet character} defined by the following properties:
\begin{enumerate}
\item $\exists \ q \in \mathbb{N} : \chi(n) = \chi(n+q) \ \forall n$;
\item if $\operatorname{gcd}(n,q)> 1$  then $\chi(n) = 0$, otherwise $\chi(n) \neq 0$;
\item $\chi(mn) = \chi(m) \chi(n) \ \forall n,m \in \mathbb{Z} $.
\end{enumerate}
For a certain $\chi(n)$, $q$ is said to be the modulus of $\chi$ and for a fixed $q$ there are several characters labelled $j$. A character is called principal and labelled $j=1$ if

\begin{equation}\label{principal}
\chi(n) = 
\left\{ 
  \begin{array}{l l}
   1 \qquad \mbox{if n and q are coprimes},  \\
   0 \qquad \mbox{otherwise}.  \\

  \end{array} \right.
\end{equation}
the series \eqref{DirichletL} modulo $q$ can be extended to a function defined over the whole complex plane (except for the case of principal characters, which have a pole at $s=1$) and satisfying the functional equation \cite{1} (p. 175)

\begin{subequations}
\begin{equation}
\Lambda(1-s, \overline{\chi}) = \frac{i^a \sqrt{q}}{\tau(\chi)} \Lambda(s,\chi),
\end{equation}
\begin{equation}
\Lambda(s,\chi) = \left( \frac{q}{\pi} \right)^\frac{s+a}{2} \Gamma{\left( \frac{s+a}{2} \right)} L(s,\chi),
\end{equation}
\begin{equation}
a = \tfrac1{2}(1-\chi(-1)),
\end{equation}
\begin{equation}
\tau(\chi) = \sum_{n=1}^q {\chi{(n)} e^{2 \pi i n/q}}.
\end{equation}
\end{subequations}

The generalized Riemann hypothesis (GRH) states that

\begin{center}
\emph{there are no Dirichlet L-function with zeros $\operatorname{Re}(\rho)>1/2$}.
\end{center}
The ordinary RH is included because $\zeta(s)$ is the Dirichlet L-function with modulus 1.

It is a well-known result that $L(s,\chi_1)$ can be expressed in terms of the zeta function, so they depend on the ordinary RH: assume $\chi_1$ is a primitive character modulo $q$ so that \eqref{principal} yields

\begin{equation*}
L(s,\chi_1) = \zeta(s) - q^{-s} \zeta(s) -\sum_P{\sum_{n=1}^{\infty}{\frac1{(P-q+nq)^s}}},
\end{equation*}
where $0<P<q$ are not coprimes with $q$. We define $p=\operatorname{gcd}(q-P,q) \in \mathbb{P}$ to obtain
\begin{equation*}
L(s,\chi_1) = q^{-s} \left[ q^{s} \zeta(s) - \zeta(s) -\sum_P{p^s\sum_{n=1}^{\infty}{\frac1{(np+p\frac{P-q}{q})^s}}} \right].
\end{equation*}

The point here is to note that $p\frac{P-q}{q}$ runs from $1$ to $p-1$, so terms $np$ are missing in the series. To complete the zeta series, we call $|p|$ the number of distinct nontrivial prime factors in the decomposition of $q$ and write

$$ -\zeta(s) = (|p|-1)\zeta(s) - |p|\zeta(s) $$

to finally obtain

\begin{equation}
L(s,\chi_1) = q^{-s} \left[ q^{s} + |p|-1 - \sum_p{p^s} \right] \zeta(s).
\end{equation}
These functions vanishes for certain periodic values $s_n = it_n$, but the nontrivial zeros are those of the zeta function.

The Hadamard product theorem states that any entire function $f(s)$ can be expressed as a product over its zeros $\rho$; for entire functions of order 1 \cite{2}:
\begin{equation}\label{generalproduct}
f(s) = s^{m_0} e^{A+Bs} \prod_{\rho}{\left( 1-\frac{s}{\rho} \right)e^{s/\rho}},
\end{equation}

where $m_0$ is the multiplicity of $\rho_0 = 0$. To apply this formula to Dirichlet L-functions $j \neq 1$, we henceforth assume $L(0,\chi) \neq 0$ (the functional equation excludes the possibility $L(0, \chi)=0$ for primitive characters because Dirichlet proved $L(1, \chi) \neq 0$.)

\begin{equation*}
\Lambda(s,\chi) = e^{A_\chi+B_\chi s} \prod_{\rho_\chi}{\left( 1-\frac{s}{\rho_\chi} \right)e^{s/\rho_\chi}}.
\end{equation*}

The logarithmic derivative gives
\begin{equation*}
\frac{\Lambda'}{\Lambda}(s,\chi) = B_\chi + \sum_{\rho_\chi}{\left( \frac1{s - \rho_\chi} + \frac1{\rho_\chi} \right)};
\end{equation*}

from the functional equation:
\begin{equation*}
\frac{\Lambda'}{\Lambda}(s,\chi) = - \overline{\frac{\Lambda'}{\Lambda}(1-\overline{s},\chi)},
\end{equation*}

where we have used $L(\overline{s},\overline{\chi}) = \overline{L(s,\chi)}$, and for $s=0$
\begin{equation*}
B_\chi = - \frac1{2} \sum_{\rho_\chi}{\left( \frac1{1 - \rho_\chi} + \frac1{\rho_\chi} \right)} = - \mathcal{T}.
\end{equation*}

Now, for the functional equation we write $\tau(\chi) = \sqrt{q} e^{i\theta}$ and $i^a e^{-i \theta} = e^{-i\theta'}$ to obtain
\begin{equation*}
\Lambda(1-s, \overline{\chi}) e^{-i\theta'/2} = \Lambda(s,\chi)e^{i\theta'/2},
\end{equation*}
so that a real-valued function $F(t,\chi) = \Lambda(1/2+it,\chi)e^{-i\theta'/2}$ can be constructed:

\begin{equation*}
F(t,\chi) = e^{(2A_\chi+\mathcal{T}-i\theta')/2} B(t,\chi) \in \mathbb{R}.
\end{equation*}

$t=0$ gives $F(0,\chi) = e^{(2A_\chi+\mathcal{T}-i\theta')/2} \in \mathbb{R}$ and
\begin{equation*}
B(t,\chi) \in \mathbb{R} \Longrightarrow \mbox{GRH}.
\end{equation*}

The GRH has profound implications; for instance, it automatically proves Goldbach's weak conjecture \cite{3}.

\section{Further extensions of the RH}

The extended Riemann hypothesis (ERH) is related to the Dedekind zeta functions

\begin{equation}\label{Dedekindzeta}
\zeta_K(s) = \sum_{I \subseteq \mathcal{O}_K} {\frac1{(N_{K/Q}(I))^s}},
\end{equation}

where $K$ is an algebraic number field, $I$ ranges through the non-zero ideals of the ring of integers $\mathcal{O}_K$ and $N_{K/Q}(I)$ is the absolute norm of $I$.

It satisfies a functional equation of the form $\Lambda_K(s) = \Lambda_K(1-s)$ \cite{1} (p. 217), so a real-valued function can be defined containing the nontrivial zeros of these functions; if a Hadamard product representation of the form \eqref{zetaHadamard} can be found, would the ERH be proven? Unfortunately, not all functions satisfying both requirements satisfy the RH: consider the Epstein zeta functions for positive definite quadratic forms:
\begin{equation}\label{Epstein}
Z(s) = \sum_{(n,m) \neq (0,0)}{\frac1{(an^2+bnm+cm^2)^s}}.
\end{equation}

The corresponding RH is known to fail for Epstein zeta functions with class number $\geq 2$, even satisfying a functional equation \cite{4}: 
\begin{subequations}
\begin{equation}
\Lambda(s) = \Lambda(1-s),
\end{equation}
\begin{equation}\label{Epsteinfunctional}
\Lambda(s) = s(s-1)\left( \frac{\sqrt{\Delta}}{2\pi} \right)^s \Gamma{(s)}Z(s),
\end{equation}
\begin{equation}
\Delta = 4ac-b^2.
\end{equation}
\end{subequations}

The problem is that these functions have zeros for $\operatorname{Re}(s)>1$, as can be seen from \eqref{Epstein}, so our proof for the original RH does not necessarily hold. This is not the case for the Dedekind zeta functions, so we expect the ERH to still hold. We also know that these Epstein functions do not have an Euler product, as it does not allow zeros for $\operatorname{Re}(s)>1$; it is widely believed that a certain function satisfies the analogous of the RH if and only if it has an Euler product, and this might be the reason.

\end{document}